\documentclass[a4paper,UKenglish,cleveref, autoref, thm-restate]{lipics-v2021}

\pdfoutput=1 

\title{One \texorpdfstring{$n$}{n} Remains to Settle the Tree Conjecture} 

\author{Jack Dippel}{McGill University, Montreal QC H3A 0G4, Canada, \and \url{https://www.cs.mcgill.ca/~jdippe/}}{jack.dippel@mail.mcgill.ca}{https://orcid.org/0000-0002-8087-3009}{}
\author{Adrian Vetta}{McGill University, Montreal QC H3A 0G4, Canada, \and \url{https://www.math.mcgill.ca/vetta/}}{adrian.vetta@mcgill.ca}{}{}

\authorrunning{J. Dippel and A. Vetta}

\Copyright{Jack Dippel and Adrian Vetta CC-BY} 



\nolinenumbers 

\usepackage{tikz}
\usepackage{caption}
\usepackage{subcaption}
\usepackage{todonotes}
\allowdisplaybreaks

\graphicspath{{images/}}

\usepackage{pdfpages}

\usepackage{subfiles} 

\begin{document}

\maketitle
	
\begin{abstract}
In the famous {\em network creation game} of Fabrikant et al.~\cite{FLM03} a set of agents play a game to build a connected graph. 
The $n$ agents form the vertex set $V$ of the graph and each vertex
$v\in V$ buys a set $E_v$ of edges inducing a graph $G=(V,\bigcup\limits_{v\in V} E_v)$. The private objective of each vertex is
to minimize the sum of its building cost (the cost of the edges it buys) plus its connection cost (the total distance from itself to every other vertex). Given a cost of $\alpha$ for each individual edge, a long-standing conjecture, called the {\em tree conjecture}, states that if $\alpha > n$ then every Nash equilibrium graph in the game is a spanning tree. After a plethora of work, it is known that the conjecture holds for any $\alpha>3n-3$. In this paper we prove the tree conjecture holds for 
$\alpha>2n$. This reduces by half the open range for $\alpha$ with only $[n, 2n)$ remaining in order to settle the conjecture.
\end{abstract}

\newpage

\section{Introduction}\label{sec:intro}

A foundational motivation for the field of algorithmic game theory was to understand the evolution and functionality of networks, specifically, the internet; see Papadimitriou~\cite{Pap01}.
Of particular fascination concerned how the actions of self-motivated agents affected the structure of the world wide web
and social networks more generally. An early attempt to study this conundrum was undertaken by Fabrikant et al.~\cite{FLM03} with their now classical {\em network creation game}.
Despite its extreme simplicity, their model (detailed below) has become highly influential for two reasons.
First, it inspired the development of a wide range of network formation models.
Second, it has lead to one of the longest-standing open problems in
algorithmic game theory, namely, the tree conjecture.
The latter motivates this research. So let's begin by describing the model and the conjecture.

\subsection{The Network Creation Game}
Consider a set of vertices, $V=\{1,2,\dots, n\}$, who attempt to construct
a connected graph between themselves. To do this, each vertex (agent) can
purchase individual edges for a fixed cost of $\alpha$ each.
Consequently, a strategy for vertex $v\in V$ is a set of (incident) edges $E_v$. Together the strategies of the agents forms a graph $G = (V,E)$, where $E=E_1\cup E_2\cup \cdots \cup E_n$.
In the {\em network creation game}, the objective of each vertex is to minimize the sum of its building and its connection cost. The {\em building cost} for vertex $v$ is $\alpha\cdot |E_v|$, the cost of all the edges it buys. The {\em connection cost} is $D(v)=\sum_{u: u\neq v} d_G(u,v)$, the sum of the distances in $G$ of $v$ to every other vertex, where $d_G(u,v)=\infty$ if there is no path between $u$ and $v$. 
That is, the total cost to the vertex is
$c_v(E) = \alpha\cdot |E_v| + D(v)$.

Given the different objectives of the agents, we study Nash Equilibria (NE) in the network creation game. A Nash equilibrium graph is a graph $G=(V,E)$ in which no vertex $v$ can reduce its total cost by changing its strategy, that is, by altering the set of edges it personally buys, given
the strategies of the other vertices remain fixed.
Thus $E_v$ is a best response to $(E_u)_{u\neq v}$, for every vertex $v$.
Observe that every Nash equilibrium graph must be a connected graph. Attention in the literature has focused on 
whether or not every Nash equilibrium graph is minimal, that is, a spanning tree.

\subsection{The Tree Conjecture}

The network creation game was designed by Fabrikant et al.~\cite{FLM03}. They proved the {\em price of anarchy} for the game is $O(\sqrt{\alpha})$, and conjectured that it is a constant. This conjecture originated from their proof that any network equilibrium graph which forms a tree costs at most $5$ times that of a star, i.e. the optimal network. They proposed that for $\alpha$ greater than some constant, every Nash equilibrium graph is a spanning tree. This was the original {\em tree conjecture} for network creation games.

In the subsequent twenty years, incremental progress has been made in determining the exact range of $\alpha$ for which all Nash equilibria are trees. Albers et al. \cite{AEE14} demonstrated that the conjecture holds for $\alpha \geq 12n \log n$. However, they also provided an counterexample to the original tree conjecture. Moreover, Mamageishvili et al.~\cite{MMM15} proved the conjecture is false for $\alpha \leq n-3$.

As a result, a revised conjecture took on the mantle of the {\em tree conjecture}, namely that every Nash equilibrium is a tree if $\alpha > n$.
Mihal\'ak and Schlegel~\cite{MS13} were the first to show that this tree conjecture holds for  $\alpha\ge c\dot n$ for a large enough constant
$c$, specifically, $c>273$. Since then, the constant has been improved repeatedly, by Mamageishvili et al.~\cite{MMM15}, then \`Alvarez and Messegu\'e~\cite{AM17}, followed by Bil\`o and Lenzner~\cite{BL19}, and finally 
by Dippel and Vetta \cite{DV22} who proved the result for $\alpha > 3n-3$.

We remark that extensions and variations of the network creation game have also been studied; we refer the interested reader to~\cite{ADH13,BGL16,BGL14,CLM17,CL15,DHM12,Q22,AM23}.

\subsection{Our Contribution}
In this paper we improve the range in which the tree conjecture is known to hold from $\alpha > 3n-3$ \cite{DV22} to $\alpha > 2n$:
\begin{theorem}\label{main}

If $G$ is a Nash equilibrium graph for the network creation game $(n,\alpha)$ and $\alpha > 2n$, $G$ is a tree.
\end{theorem}

As alluded to in the title,
this result decreases the length of the open range for $\alpha$
to just $n$.

Our high-level strategy to prove Theorem~\ref{main} is straightforward: we assume the existence a of biconnected component $H$ in a Nash equilibrium and then prove, via consideration of a collection of strategy deviations, that some vertex has a better strategy, provided $\alpha > 2n$. This contradicts
the best response conditions and proves that every Nash equilibrium graph is a spanning tree.

Our approach applies some prior methodologies and combines them with some original tools. One important technique we exploit is that of {\em min-cycles}, introduced by 
Lenzner~\cite{L14}. A min-cycle is a cycle in a graph
with the property that, for every pair of vertices in the cycle, the cycle contains a shortest path between the pair. For certain values of $\alpha$, currently $\alpha> 2n-3$, min-cycles are known to have the nice property that each vertex in the cycle buys exactly one edge of the cycle. This property has been leveraged in various ways~\cite{AM17,BL19,DV22} to show that no min-cycles can exist in a Nash equilibrium graph for large enough $\alpha$. Importantly, the smallest cycle in any graph must be a min-cycle. 
Consequently, for large enough $\alpha$, 
there can be no smallest cycle and, thus, no cycle at all.
This implies that, for such an $\alpha$, all Nash equilibrium graphs are trees.

The techniques we develop concern the analysis of the
presupposed biconnected component $H$ in the Nash equilibrium. The existence of $H$ implies the existence of a special vertex, called $r$, which has the lowest connection cost amongst all vertices in $H$. That is, $D(r) \leq D(v)$ for all $v \in H$. Given $r$, we take a shortest path tree $T$ in $G$ rooted at $r$. 
The key is to exploit the structural properties inherent in $T$. To derive these properties, we design a new class of strategy deviations available to vertices in $H$. Because these deviations must be non-improving responses, they impose numerous beneficial restrictions on the edges in $H$. In particular, we can then show that every edge in $H\setminus T$ increases the sum of degrees by more than 2. This absurdity proves the Nash equilibrium graph is a tree.

\subsection{Overview of Paper}

This paper has three main sections.
Section~\ref{sec:prelim} consists of preliminaries
and contains three things. First, we examine the structure of Nash equilibrium graphs which contain a
hypothesized biconnected component $H$. Second, we discuss min-cycles, which appear in several previous papers, each of which contains a useful lemma that we take advantage of in this paper. Third, we introduce some new lemmas pertaining to $T$, the shortest path tree mentioned above.

Section~\ref{sec:strategies} has two parts. The first  presents a set of three deviation strategies. That is, three specific ways a vertex in a biconnected component might consider changing which edges it buys in order to decrease its personal cost. For a Nash equilibrium graph these alternate strategies cannot reduce the personal cost, and this allows us to derive deviation bounds for that vertex and the edges it buys. The second part uses these deviations bounds to prove claims about the structure of the graph. In particular, we show that vertices which buy edges outside $T$ cannot be too close to one another.

Finally, Section~\ref{sec:result} presents two main lemmas. The first shows that, in a set $U$ of three vertices adjacent in $T$, all at different depths, where the lowest buys an edge outside $T$, only one vertex in $U$ can have degree two in $H$. The second lemma bounds the number of times such sets can intersect. These two lemmas are the core of the proof of~\Cref{main}. After proving the main result, the paper culminates with concluding remarks on the future of the tree conjecture.

\section{Preliminaries}\label{sec:prelim}

Recall our basic approach is to prove by contradiction the non-existence of a biconnected component
in a Nash equilibrium graph. Accordingly, we begin in Section~\ref{sec:biconnected} by introducing
the notation necessary to analyse biconnected components in network creation game graphs.
To disprove the existence of a biconnected component it suffices to disprove the existence of 
a cycle. Of particular importance here is the concept of a min-cycle. So in Section~\ref{sec:min-cycle}
we present a review of min-cycles along with a corresponding set of fundamental lemmas which appear in~\cite{L14,BL19,MMM15}. Finally, in Section~\ref{sec:technical}, we prove a
collection of technical results concerning the shortest path tree $T$
that we will utilize throughout the rest of the paper. 
\Cref{deg2} is an especially useful tool for this method, and could prove valuable in future work on this problem.

\subsection{Biconnected Components}\label{sec:biconnected}

Given a subgraph $W$ of a graph $G$, we let $d_W(v,u)$ denote the distance between $v$ and $u$ in $W$. 
In particular, $d(u,v)=d_G(u,v)$ is the distance from $v$ to $u$ in the whole graph. Let 
$D(v)=\textstyle\sum_{u: u\neq v} d_G(u,v)$ be the {\em connection cost} for vertex $v$,
the sum of the distances from $v$ to every other vertex.

In a Nash equilibrium graph $G$ containing biconnected components, we will refer to the largest biconnected component as $H$. For a vertex $v \in H$, we denote by $C(v)$ a smallest cycle containing $v$, breaking ties arbitrarily. Let $r\in H$ be a vertex whose connection cost is smallest amongst all the vertices in $H$. 
Once built, an edge $uv$ of $G$ can be traversed in either direction. Thus $G$ is an undirected graph.
However, it will often be useful to view $G$ as a directed graph. Specifically,
we may orient $uv$ is from $u$ to $v$ if the edge was bought by $u$
and orient it from $v$ to $u$ if the it was bought by $v$.

Given $r$, we will make heavy use of the shortest path tree $T$ rooted at $r$. 
There may be multiple options for the choice of $T$. So we insist our choice of $T$ has the following property: For all shortest paths $P$ between $r$ and any vertex $u$ which is directed from $r$ to $u$, we have $P \subseteq T$. The proof of \Cref{directedcycles} shows that we will never have two shortest paths directed from $r$ to the same vertex $v$, so this choice of $T$ is well-defined. 
Next, given an edge $uv$ where $u$ is the parent of $v$ in $T$, we say that 
$uv$ is a {\em down-edge} of $T$ if the edge was bought by $u$; otherwise we say it is an {\em up-edge} of $T$.
We denote by $T^{\downarrow}$ the set of down-edges and by $T^{\uparrow}$ the set of up-edges of $T$.
Let $T(v)$ be the subtree rooted at $v$ in the tree $T$ rooted at $r$.
If $uv\in T^{\downarrow}$ then we define $T(uv) = T(v)$. Otherwise, if $uv\notin T^{\downarrow}$ then we define $T(uv) = \emptyset$.

Finally, there remain two types of sets we need to define. The first type is the ``$S$-sets", which are sets of vertices. In a graph $G$, with a subgraph $W$, $S_W(v)$ is the set of vertices with a shortest path in $G$ to any vertex in $W$ which contains $v$. Therefore $S_W(v)$ will always contain $v$, and it will contain no vertices of $W\setminus v$. $W$ can even be a single vertex $u$, i.e. $S_u(v)$ is the set of all vertices in $G$ with a shortest path to $u$ containing $v$. 
The second type is the ``$X$-sets", which are sets of edges. For largest biconnected component $H$ of $G$
and the shortest path tree $T$, the set $X_0$ is the set of all edges in $H\setminus T$, or {\em out-edges}. For all integers $i \geq 1$, $X_i$ contains the set $X_{i-1}$ as well as all edges $uv \in H\cap T^{\downarrow}$ bought by the parent $u$ and where the child $v$ buys an edge of $X_{i-1}$. The $X_i$ sets can be thought of as the set of all out-edges and all down-edges in directed paths of length $\leq i$ to a vertex which buys an out-edge. Sometimes it will be useful to include in the $X$-sets all the edges $vu \in H\cap T^{\uparrow}$, bought by child $v$ to parent $u$. To do this, we add superscript $+$ to the $X$ set. For instance $X_0^+=X_0\cup T^{\uparrow}$ is the set of all up-edges and out-edges, which happens to be the set of all edges $e$ with $T(e)=\emptyset$. Similarly, $X_i^+=X_i\cup T^{\uparrow}$.

\begin{figure}
\centering
\includegraphics[]{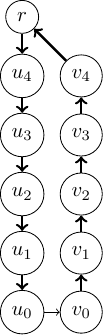}
\caption{In the diagram, edges of $T$ are bold. Vertex $u_0$ buys an edge $u_0v_0 \in X_j$ for all $j \geq 0$. Similarly, each vertex $u_i$ buys edge $u_iu_{i-1} \in X_j$ for all $4 \geq i \geq 1$, $j \geq i$. Each vertex $v_i$ buys $v_iv_{i+1} \in X_j^+$ for all $3 \geq i \geq 0$, $j \geq 0$. $v_4$ buys $v_4r \in X_j^+$ and $r$ buys $ru_4 \in X_k$ for $j \geq 0$, $k \geq 5$.}
\label{fig:xsets}
\end{figure}

\subsection{Min-Cycles}\label{sec:min-cycle}

Recall, a {\em min-cycle} is a cycle in a graph with the property that, for every pair of vertices in the cycle, the cycle contains a shortest path between the pair.
This concept, introduced by Lenzner~\cite{L14}, has proved very fruitful in the study of 
Nash equilibria in network creation games. In particular, we now present three useful
lemmas concerning min-cycles that appear in~\cite{L14,BL19,MMM15}, respectively.
For completeness and because some of the ideas used are informative, we include short proofs of these three lemmas; similar proofs appeared in~\cite{DV22}

\begin{lemma}\label{mincycleedges}\cite{L14}	
The smallest cycle $C$ containing an edge $e$ is a min-cycle.
\end{lemma}

\begin{proof}
Consider the smallest cycle $C$ containing an edge $e$. Suppose for the sake of contradiction that there are two vertices $u,v\in C$ such $d_G(u,v) < d_C(u,v)$. Without loss of generality, suppose the shortest path between $u$ and $v$, labelled $P$, lies entirely outside $C$. Note that $C$ contains two paths
between $u$ and $v$. Let $Q$ be the path of $C$ from $u$ to $v$ that contains $e$. 
Then $P \cup Q$ is a cycle containing~$e$ that is strictly smaller than $C$, a contradiction.

\end{proof}

A nice property of min-cycles is that every vertex in the min-cycle buys
exactly one edge of the cycle.
\begin{lemma}\label{directedcycles}\cite{BL19}
Let $\alpha > 2(n-1)$. Every min-cycle in a Nash equilibrium graph $G$ is directed.
\end{lemma}

\begin{proof}
Let $C$ be a min-cycle that is {\em not} directed. Then there is a vertex $v$ that buys two edges of the cycle, say $vx$ and $vy$. Let $uw \in C$ be an edge furthest from $v$ (if $|C|$ is odd, there is a unique choice, otherwise choose on of two furthest edges). Let $u$ buy $uw$. We've chosen $uw$ so that both $u$ and $w$ have a shortest path to $v$ which does not contain $uw$, therefore every vertex in $G$ has shortest path to $v$ without $uw$. Without loss of generality, let $d(u,y) < d(u,v)$. Then because of how we chose $uw$, both $v$ and $x$ have a shortest path to $u$ which does not contain $vx$, therefore every vertex in $G$ has shortest path to $u$ without $vx$.

First, consider that if $u$ sells $uw$ and buys $uv$, no vertices become farther from $v$. 
It follows, by the Nash equilibrium conditions, that 
\begin{equation}\label{eq:directed-cycle-1}
D(u) \leq D(v) + n-1
\end{equation}
because $D(u)$ must be less than the cost to swap $uw$ for $uv$ and use the edge $uv$ in followed by the shortest path from $v$ to each of the $n-1$ vertices. 

On the other hand, $v$ can sell {\em both} $vx$ and $vy$ and instead buy the edge $vu$ without increasing the distance from $u$ to any other vertex. All of the vertices which used $vy$ in their shortest paths to $u$ become closer to $u$.  
It follows, by the Nash equilibrium conditions, that 
\begin{equation}\label{eq:directed-cycle-2}
D(v) \leq D(u) + n-1 - \alpha
\end{equation}

Together (\ref{eq:directed-cycle-1}) and  (\ref{eq:directed-cycle-2}) give $D(v) \leq D(v) + 2(n-1) - \alpha < D(v)$.
This contradiction implies that $C$ must be a directed cycle.

\end{proof}

Furthermore, we may now obtain a lower bound on the girth of a Nash equilibrium graph in 
terms of the cost $\alpha$.
\begin{lemma}\label{mincyclesize}\cite{MMM15}
Any cycle $C$ in a Nash equilibrium graph has $|C| \geq \frac{2\alpha}{n}+2$.
\end{lemma}

\begin{proof}
	Take a minimum length cycle $C = v_0, v_1\dots, v_{k-1}, v_{k}=v_0$ in $G$. 
	By \Cref{mincycleedges}, $C$ is a min-cycle. Therefore, by \Cref{directedcycles}, $C$ is a directed cycle.
	So we may assume $e_i=v_iv_{i+1}$ is bought by $v_i$, for each $0 \leq i \leq k-1$. 
	Now, for each vertex $u\in V$, we define a set $L_u\subseteq \{0,1,\dots, k-1\}$ as follows. We have $i\in L_u$ if and only if
	{\em every} shortest path from $v_i\in C$ to $u$ uses the edge $e_i$.
	
	We claim $|L_u| \le \frac{|C|-1}{2}$ for every vertex $u$.
	If not, take a vertex $u$ with $|L_u| > \frac{|C|-1}{2}$. 
	Let $d(v_i,u)$ be the shortest distance between $u$ and $v_i\in C$. 
	Next give $v_i$ a label $\ell_i= d(v_i,u)-d(v_{i+1},u)$. Observe that $\ell_i\in \{-1,0,1\}$.
	Furthermore, the labels sum to zero as
	$$\textstyle\sum_{i=0}^{|C|-1} \ell_i \ =\  \textstyle\sum_{i=0}^{|C|-1} \left( d(v_i,u)-d(v_{i+1},u) \right) 
		\ =\ \textstyle\sum_{i=0}^{|C|-1} d(v_i,u) - \textstyle\sum_{i=1}^{|C|} d(v_{i},u)  \ = \ 0$$
	Now take a vertex $v_i$ in $C$ that uses $e_i$ in {\em every} shortest path to $u$; that is, $i\in L_u$.
	Then $\ell_i=1$ and $\ell_{i-1}\ge 0$. 
	In particular, if $|L_u| > \frac{|C|-1}{2}$ then there are $> \frac{|C|-1}{2}$ positive 
	labels and $> 1+\frac{|C|-1}{2}$ non-negative labels.
	Hence, there are $< |C|-\left(1+\frac{|C|-1}{2}\right) = \frac{|C|-1}{2}$ negative labels.
	But then the sum of the labels is strictly positive, a contradiction.
	
	Now, as $|L_u| \le \frac{|C|-1}{2}$ for every vertex $u$, there must exist a $v_i,e_i$ pair where $v_i$ needs $e_i$ for its shortest paths to $\leq \frac{n}{2}$ vertices. For the vertices that do require $e_i$ in their shortest paths, deleting $e_i$ increases their distance by $\leq |C|-2$, as we can replace $e_i$ with $C\setminus e_i$ in those paths. Therefore, if $v_i$ sells $e_i$ its cost increases by $\leq (|C|-2)\cdot \frac{n}{2} -\alpha$. This must be non-negative by the Nash equilibrium conditions.
	Rearranging, we have $\frac{2\alpha}{n} + 2\leq |C|$. 
 
\end{proof}

An important immediate consequence of~\Cref{mincyclesize} is that we may assume that the
girth of any Nash equilibrium graph of interest here is at least $7$.
\begin{corollary}\label{7cycle}
Let $\alpha > 2n$. All cycles in a Nash equilibrium graph have length at least $7$.
\end{corollary}

\subsection{The Shortest Path Tree \texorpdfstring{$T$}{T}}\label{sec:technical}

Let's now consider the shortest path tree $T$ rooted at the vertex $r$ with smallest connection cost
amongst the vertices in the biconnected component $H$.
In this section, we present some technical lemmas that give useful insights into the structure of this tree $T$.
We begin by upper bounding the size of any subtree in in $T$.
\begin{lemma}\label{maxn2}
If $v\not=r$ then $|T(v)|\leq \frac{n}{2}$.
\end{lemma}

\begin{proof}
Recall $T$ is the shortest path tree of the vertex $r$ with the
smallest connection cost and $T(v)$ is the subtree of $T$ rooted at $v$.
For any vertex $x \in T(v)$, it immediately follows that
$d(v,x) = d(r,x)- d(v,r)$. On the other hand, for any vertex
$x \not\in T(v)$, the triangle inequality implies that
$d(v,x) \leq d(r,x) + d(v,r)$. Thus
\begin{eqnarray*}
D(v) 
&\leq& D(r) - d(v,r)\cdot |T(v)| + d(v,r)\cdot (n-|T(v)|) \\
&=& D(r) + d(v,r)\cdot (n-2\cdot|T(v)|)
\end{eqnarray*}
But $D(r) \leq D(v)$. Therefore $d(v,r)\cdot (n-2\cdot|T(v)|)\geq 0$ and rearranging gives $|T(v)| \leq \frac{n}{2}$.

\end{proof}

The next two lemmas concern the properties of paths related to
any vertex $u$ that buys an edge $uv \in X_i$ where $i \leq 2$.
\begin{lemma}\label{altpath}
Let $\alpha > 2n$. If $u$ buys $uv \in X_i$, for some $i \leq 2$, then all $w \in T(uv)$ have a path to $r$ of length $\leq d(w,r)+2i$ which does not contain $u$. 
\end{lemma}

\begin{proof}
  If $uv \in X_0$ then, as $uv\notin T^{\downarrow}$, we have $|T(uv)|=\emptyset$. Thus the claim is trivially true. 

   So assume $uv \not\in X_0$. Then $uv$ is the first edge in a path of length $i$ in $T$ from $u$ to a vertex $x$ which buys $xy \in H \setminus T$. Thus, there is a path $P_0$ in $T$ of length $\ell(P_0)=i-1$ from $v$ to $x$. 
   Now let $P_1$ be the path from $y$ to $r$ in $T$. Observe that $P_1$ does not contain $u$; otherwise $xy$ and $T$ would define a cycle of length $\leq 6$ (because $i\le 2$), contradicting \Cref{7cycle}. Note also that $P_1$ has length $\ell(P_1)\leq d(r,v)+i$ because $T$ is a shortest path tree.

    Next, let $P_2$ be the path from $w$ to $v$. Note $P_2$ has length $d(w,r) - d(v,r)$ and also does not contain $u$. Consequently $P_2,P_0,xy,P_1$ is a path from $w$ to $r$. This path has length at most
    \begin{eqnarray*}
    \ell(P_2)+\ell(P_0)+ 1+\ell(P_1) &\le&
    \left(d(w,r) - d(v,r)\right) + \left(i-1\right)+ 1 + \left( d(v,r)+i\right) \\
    &=& d(w,r) + 2i
       \end{eqnarray*}
      Furthermore this path does not contain $u$. 
    
\end{proof}

\begin{lemma}\label{X2position}
Let $\alpha > 2n$. If $u$ buys $uv \in X_i$, for some $i \leq 2$, then $d(u,r) \geq \lfloor\frac{|C(u)|}{2}\rfloor-i$.
\end{lemma}

\begin{proof}
So $uv \in X_i$ is associated with an edge $xy\in X_0$, where $xy=uv$ if $uv \in X_0$. Let $C$ be the cycle defined by $T+xy$. Because $T$ is a shortest path tree it cannot be the case that $y$ is an ancestor of $x$ in $T$. Neither is $y\in T(u)$; if so, $d(u,y)\le 3$ and therefore $\ell(C) \leq 6$, contradicting \Cref{7cycle}. Thus $y \not \in T(u)$ and, consequently, $uv \in C$.

Let $z$ be the lowest common ancestor of $x$ and $y$.
Then $C$ is a path from $z$ to $x$, plus the edge $xy$, plus a path from $y$ to $z$. So it has length 
$$d(x,z) + 1 + d(y,z) \le 2d(x,r)+2 = 2d(u,r)+2i+2$$ 
Observe that to achieve this maximum length, we must have $d(x,z) = d(x,r) = d(y,r)-1$, and $C$ must consist of two shortest paths from $y$ to $r$, meaning $C$ will also contain~$r$. 

Next consider the smallest cycle $C(u)$ containing $u$. Suppose for the sake of contradiction that $d(u,r) \leq \lfloor\frac{|C(u)|}{2}\rfloor -i-1$. If $|C(u)|$ odd, then $2d(u,r)+ 2i + 3 \leq |C(u)|$, which contradicts the choice of $C(u)$. Therefore, $|C(u)|$ is even and $2d(u,r) + 2 + 2i \leq |C(u)|$. However, $|C| \leq 2d(u,r)+2i +2$.
Thus $|C| = |C(u)|$ and therefore $C$ a min cycle, by \Cref{mincycleedges}. Furthermore, $C$ is directed by \Cref{directedcycles}. This is a contradiction because then $xy$ is in a directed shortest path, and must be part of $T$ and not an edge in $X_0$. Hence $d(u,r) \geq \lfloor\frac{|C(u)|}{2}\rfloor-i$, as desired.

\end{proof}

We present one last technical lemma in this section. We remark that this lemma is of critical value 
in our analysis and, we believe, may be of importance in achieving future improvements.

\begin{lemma}\label{deg2}
        If $uvw$ is a directed path in $H$ and $\deg_H(v) =2$, then $|S_H(v)| \geq \frac{\alpha}{2(|C(v)|-3)}$. If, in addition, $uvw$ is a directed path of down-edges in $T$, then $|S_H(v)| \geq |T(w)|$.
\end{lemma}

\begin{proof}
    Consider the strategy change of $u$ selling $uv$ and buying $uw$. As $\deg_H(v)=2$, the only vertices which become farther from $u$ are those of $S_H(v)$. They all become farther from $u$ by 1. All the vertices in $S_u(w)$ conversely, become closer to $u$ by 1. Therefore, by the equilibrium conditions, we must have $|S_H(v)| \geq |S_u(w)|$, otherwise $u$ has an incentive to switch strategies. Further, if $uvw$ is a path of down-edges in $T$, then because $u$ is the ancestor of $w$ and $T$ is a shortest path tree,  $T(w) \subseteq S_u(w)$. Thus $|S_H(v)| \geq |T(w)|$, as desired.

    Now consider a second strategy change in which $u$ deletes $uv$. There are two sets of vertices whose distance to $u$ increases. Shortest paths from $u$ to $S_H(v)$ go from using $uv$ to using $C(v)\setminus uv$; note that as $\deg_H(v) =2$ it must be that case that $uv, vw\in C(v)$.
    The distance of these vertices from $u$ then increase by $|C(v)|-2$. Meanwhile, paths from $u$ to $S_u(w)$ go from using $uvw$ to using at most $C(v)\setminus\{uv,vw\}$, an increase in distance of at most $ |C(v)|-4$. No other vertices have increased distance to $u$. Therefore $u$'s cost as a result changes by at most
 \begin{align*}
    (|C(v)|-2)\cdot |S_H(v)| + (|C(v)|-4)\cdot |S_u(w)|-\alpha
    \leq (2|C(v)|-6)\cdot |S_H(v)| -\alpha 
\end{align*}
Here the inequality holds as $|S_H(v)| \geq |S_u(w)|$.
This change must non-negative by the equilibrium condition. This implies $|S_H(v)| \geq \frac{\alpha}{2(|C(v)|-3)}$, as desired.

\end{proof}

\section{A Class of Deviation Strategies}\label{sec:strategies}

We now present a new class of deviation strategies. Specifically, in Section~\ref{sec:three-strategies} 
we study three related strategies that a vertex in the biconnected component may ponder. 
Concretely, we derive bounds on the change in cost to the vertex resulting from using these strategy 
changes. In particular, we will use the fact that, at a Nash equilibrium, these cost changes must
be non-negative.
Then, in Section~\ref{sec:observations}, we apply these deviation strategies to make three substantial 
observations concerning vertices that buy multiple edges of specific types in $H$.

\subsection{Three Deviation Strategies}\label{sec:three-strategies}

\begin{figure}
\centering
\begin{subfigure}{0.41\textwidth}
    \scalebox{0.55}{\includegraphics[]{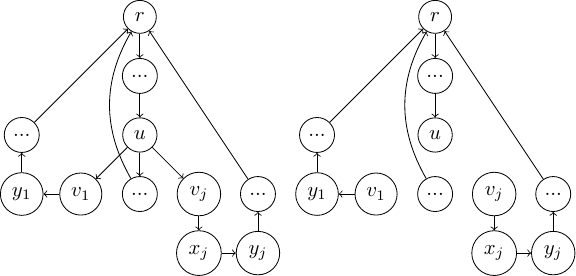}}
    \caption{Application of \Cref{mainstrategy1}. $u$ sells its edges to $v_1$,..., $v_j$}
    \label{fig:str1}
\end{subfigure}
\hfill
\begin{subfigure}{0.41\linewidth}
         \centering
         \scalebox{0.55}{\includegraphics[]{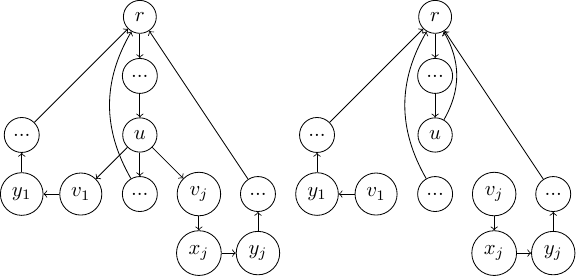}}

\caption{Application of \Cref{mainstrategy2}. $u$ sells its edges to $v_1$,..., $v_j$ and buys edge $ur$.}
    \label{fig:str2}
\end{subfigure}
\hfill
\begin{subfigure}{0.41\linewidth}
         \centering
         \scalebox{0.55}{\includegraphics[]{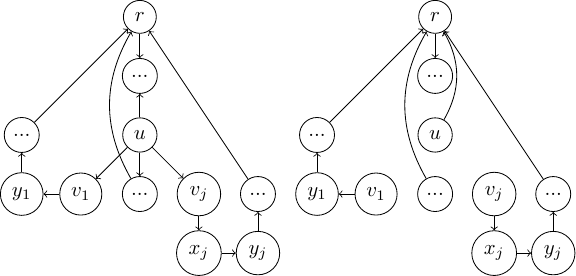}}
    \caption{Application of \Cref{mainstrategy3}. $u$ sells its edges to $v_1$,..., $v_j$ and buys edge $ur$.}
    \label{fig:str3}
    \end{subfigure}
        
\caption{The Three Main Strategy Changes for Vertex $u$}
\label{fig:figures}
\end{figure}

As stated, we now introduce a class of deviation strategies. The first, shown in \Cref{fig:str1}, involves a vertex $u$ selling edges of $X_2$. \Cref{altpath} guarantees that the graph is still connected without these edges. Now $u$'s cost will change: it saves $\alpha$ for each edge sold, but its distance to many vertices may increase by varying amounts, requiring us to bound the new total distance to all vertices. The second deviation, depicted in \Cref{fig:str2}, is very similar to the first. The only change is that $u$ buys the edge $ur$. This reduces the saving by $\alpha$ for the extra edge bought, but it also reduces the bound on the increased distance to the other vertices considerably. Finally, the third deviation, pictured in \Cref{fig:str3}, is slightly different than the second. The bound given is weaker, because $u$ may sell edges of $X_2^+$, rather than just $X_2$. Selling the edge to $u$'s parent reduces some of the savings from the previous strategies.

\begin{lemma}\label{mainstrategy1}
If $\alpha > 2n$, $u$ buys $uv_1 \in X_{i_1},...,uv_j\in X_{i_j}$ for $j\geq 1$, $i_k \leq 2$, for all $k, 1 \leq k \leq j$, $u=u_0,u_1,...u_{d(u,r)-1},u_{d(u,r)}=r$ is the path from $u$ to $r$ in $T$, then the change of $u$ selling $uv_1,...,uv_j$ results in a change of cost of at most 
\begin{equation}\label{mse1}
    d(u,r)\cdot n - \textstyle\sum\limits_{l=0}^{d(u,r)-1}2\cdot|T(u_l)|-j \cdot\alpha+ \textstyle\sum\limits_{k=1}^j (2i_k+2d(u,r))\cdot|T(uv_{i_k})| 
\end{equation}
\end{lemma}

\begin{proof}
    We commence by bounding $D(u)$ before the changes, by comparing it to $D(r)$. By the triangle inequality, $u$'s distance to any vertex $v$ is at most $d(r,u)+d(r,v)$ because there is a path from $u$ to $r$ and a path from $r$ to $v$. This is the bound we use for vertices in $T(r)\setminus T(u_{d(u,r)-1})$. For other vertices this is obviously not the shortest path. For instance, for any vertex $v \in T(u)$, $d(u,v) = d(r,v) - d(r,u)$. Furthermore, for $l$ in $1 \leq l \leq d(u,r)$, the length of $u$'s shortest path to $v \in |T(u_l)\setminus T(u_{l-1})|$ differs from $r$'s by $(d(u,r)-2l)$. Therefore we have the bound:
\begin{equation*}
    D(u) \leq D(r) -d(u,r)\cdot|T(u)| - \textstyle\sum\limits_{l=1}^{d(u,r)}(d(u,r)-2l)\cdot|T(u_l)\setminus T(u_{l-1})|
\end{equation*}

Now $D(r) \leq D(u)$. Moreover, if we rearrange, we can more concisely write this as
\begin{equation*}
    0 \leq d(u,r)\cdot n - \textstyle\sum\limits_{l=0}^{d(u,r)-1}2\cdot|T(u_l)|
\end{equation*}
    
    To complete the proof, consider what happens when $u$ sells $uv_1,...,uv_j$. First, $u$ obviously saves $j\alpha$ from it edge costs. Second, by \Cref{altpath}, we know that all  $w\in T(uv_i)$ have a path of length  $\leq 2i_k + d(w,r)$ to $r$ that does not contain $u$. This leads to an increase in distance from $u$ to $w$ of at most
    $$(d(u,r)+2i_k+d(w,r)) - (d(w,r)-d(u,r)) = 2d(u,r) + 2i_k$$  
    This yields our desired bound on the change of cost for $u$:
    \begin{equation*}
    0 \leq d(u,r)\cdot n - \textstyle\sum\limits_{l=0}^{d(u,r)-1}2\cdot|T(u_l)|-j \cdot\alpha+ \textstyle\sum\limits_{k=1}^j (2i_k+2d(u,r))\cdot|T(uv_{i_k})| 
\end{equation*}

\end{proof}

\begin{corollary}\label{mainstrategy2}
    If $\alpha > 2n$, $u$ buys $uv_1 \in X_{i_1},...,uv_j\in X_{i_j}$ for $j\geq 1$, $i_k \leq 2$, for all $k, 1 \leq k \leq j$, $u=u_0,u_1,...u_{d(u,r)-1},u_{d(u,r)}=r$ is the path from $u$ to $r$ in $T$, then the strategy change of $u$ selling $uv_1,...,uv_j$ and buying $ur$ results in a change of cost of at most  
 \begin{equation*}
 n-|T_{u_\frac{d(u,r)}{2}}| -\textstyle\sum\limits_{l < \frac{d(u,r)}{2}}2\cdot|T(u_l)|-(j-1) \cdot\alpha+ \textstyle\sum\limits_{k=1}^j (2i_k+d(u,r)+1)\cdot|T(uv_{i_k})|
\end{equation*}
where the second term is removed if $\frac{d(u,r)}{2}$ is not integral.
\end{corollary}

\begin{proof}
    This strategy change is identical that of to \Cref{mainstrategy1} except $u$ also buys $ur$. 
    This means that the new path from $u$ to $r$ is now length 1, not length $d(u,r)$. Hence we replace $d(u,r)\cdot n$ with $n$ in \Cref{mse1}. There are still vertices which are as close or closer to $u$ then to $r$. The vertices of $|T_{u_\frac{d(u,r)}{2}}|$, which only exists if $\frac{d(u,r)}{2}$ is integral, are at-least as close to $u$ as to $r$, therefore $u$ does not use the edge $ur$ on the path to any of the vertices in $|T_{u_\frac{d(u,r)}{2}}|$, thus we subtract $|T_{u_\frac{d(u,r)}{2}}|$ from \Cref{mse1}. For the remaining sets $T(u_l)$ for $l < \frac{d(u,r)}{2}$, each time $l$ decreases by 1, the vertices in the set are 1 closer to $u$ and 1 farther from~$r$. 

    Finally, the last two terms of \Cref{mse1} are affected as well. We now buy an additional edge, so we add $\alpha$, resulting in the $(j-1)\cdot\alpha$ term. Furthermore, as previously stated. The path from $u$ to $r$ now has length 1, meaning for all  $w\in T(uv_i)$ the increase in distance from $u$ to $w$ is at most
    $$(2d(u,r) + 2i_k) - (d(u,r) - 1) =1+ d(u,r) + 2i_k $$  
    This yields our desired bound on the change of cost for $u$:
 \begin{equation*}
 n-|T_{u_\frac{d(u,r)}{2}}| -\textstyle\sum\limits_{l < \frac{d(u,r)}{2}}2\cdot|T(u_l)|-(j-1) \cdot\alpha+ \textstyle\sum\limits_{k=1}^j (2i_k+d(u,r)+1)\cdot|T(uv_{i_k})| 
\end{equation*}
    
\end{proof}

\begin{corollary}\label{mainstrategy3}
    If $\alpha > 2n$, $u$ buys $uv_1 \in X^+_{i_1},...,uv_j\in X^+_{i_j}$ for $j\geq 1$, $i_k \leq 2$, for all $k, 1 \leq k \leq j$, then the strategy change of $u$ selling $uv_1,...,uv_j$ and buying $ur$ results in a change of cost of 
        \begin{multline*}
    \leq n -(j-1)\alpha -(d(u,r)+1)\cdot |T(u)| + \textstyle\sum\limits_{k=1}^j (2i_k+d(u,r)+1)\cdot |T(uv_{i_k})| 
\end{multline*}
\end{corollary}

\begin{proof}
    This strategy change is nearly identical to \Cref{mainstrategy2}, except now $u$ is potentially selling the edge to its parent in $T$.     
    We are no longer certain that $u$ is closer to its former parent than $r$ is, nor any vertices along the original path from $u$ to $r$. 
    Thus we lose the savings of all $T(u_l)$ except for $l=0$, as $T(u_0)= T(u)$.
    
\end{proof}

\subsection{Observations arising from the Deviation Strategies}\label{sec:observations}

We now apply these three deviation strategies to make three observations
concerning vertices that buy multiple edges of specific types in $H$.
In turn these observations will be instrumental in proving the main result in Section~\ref{sec:result}.
The first observation simply states that no vertex can by two edges in $X^+_1$.

\begin{observation}\label{X1}
No vertex $u\in H$ buys two edges $uv_1, uv_2\in X_1^+$.
\end{observation}

\begin{proof}
Suppose $u\in H$ buys $uv_1,uv_2\in X_1^+$. Observe $r$ cannot buy an edge in $X_1^+$, as this would imply the existence of a short cycle, contradicting~\Cref{7cycle}; therefore $u\not=r$.  Hence, \Cref{maxn2} implies $|T(u)| \leq \frac{n}{2}$. We now apply \Cref{mainstrategy3} to the strategy change where $u$ sells $uv_1, uv_2\in X_1^+$ and buys $ur$. The change in cost for $u$ is at most
\begin{eqnarray*}
\lefteqn{n - \alpha -(d(u,r)+1)\cdot |T(u)| + \textstyle\sum\limits_{k=1}^2 (2i_k+d(u,r)+1)\cdot |T(uv_{i_k})|}\\ 
 &=& n - \alpha- (d(r,u)+1)\cdot \left(|T(u)|-|T(uv_1)|- |T(uv_2)| \right) +2\left(|T(uv_1)| + |T(uv_2)|\right)\\
 &\leq& n - \alpha + 2\cdot |T(u)|\\
    &\leq& n  -\alpha+2\cdot|T(u)| \\
    &\leq& n  -\alpha +2\frac{n}{2}  \\
    &\leq& 2n -\alpha \\
    &<& 0
    \end{eqnarray*}
This contradicts the Nash equilibrium conditions. 

\end{proof}

Next we observe that no vertex can buy three edges in $X^+_2$. Furthermore any vertex that buys
two edges in $X^+_2$ must be the parent of a child at the root of a large subtree in $T$.

\begin{observation}\label{X2}
No vertex $u\in H$ buys three edges $uv_1,uv_2,uv_3\in X_2^+$, and if $u$ buys $uv_1, uv_2 \in X_2^+$, then $|T(uv_1)\cup T(uv_2)| > \frac{n}{4}$.
\end{observation}

\begin{proof}
Suppose $u\in H$ buys $uv_1, uv_2, uv_3\in X_2^+$. Again, $r$ cannot buy edges in $X_2$, as this would contradict \Cref{7cycle}; so $u\not=r$. Thus \Cref{maxn2} implies $|T(u)| \leq \frac{n}{2}$. Now we apply \Cref{mainstrategy3} to the strategy change where $v$ sells $vv_1,vv_2,vv_3\in X_2^+$ and buys $vr$. The change in cost for $v$ is:
\begin{eqnarray*}
\lefteqn{n - 2\alpha -(d(u,r)+1)\cdot |T(u)| + \textstyle\sum\limits_{k=1}^2 (2i_k+d(u,r)+1)\cdot |T(uv_{i_k})|}\\ 
 &=& n - 2\alpha- (d(r,u)+1)\cdot \left(|T(u)|-|T(uv_1)|- |T(uv_2)|- |T(uv_3)| \right) \\&~&~ +4\left(|T(uv_1)| + |T(uv_2)|+ |T(uv_3)|\right)\\
 &\leq& n - 2\alpha + 4\cdot |T(u)|\\
    &\leq& n  -2\alpha+4\cdot|T(u)| \\
    &\leq& n  -2\alpha +4\frac{n}{2}  \\
    &\leq& 3n -2\alpha \\
    &<& 0
    \end{eqnarray*}  
This contradicts the Nash equilibrium conditions. Now suppose $u\in H$ buys $uv_1, uv_2\in X_2$. We can apply \Cref{mainstrategy3} to the strategy change where $v$ sells $uv_1, uv_2\in X_2^+$ and buys $ur$. The change in cost for $u$ is at most
\begin{align*}
&n -\alpha - (d(u,r)+1)\cdot \left(|T(v)| -|T(uv_1)|-|T(uv_2)| \right) +4\cdot  \left(|T(uv_1)| + |T(uv_2)|\right)  \\
    &\leq n -\alpha +4\cdot |T(uv_1)\cup T(uv_2)| \\
    &< 4\cdot |T(uv_1)\cup T(uv_2)| -n 
    \end{align*}
This is non-negative by the Nash equilbrium conditions.
Therefore $|T(vv_1)\cup T(vv_2)| > \frac{n}{4}$, as desired.

\end{proof}

Finally, we observe some properties that follow when a vertex buys two edges in $X_2$.

\begin{observation}\label{X2depth}
If $u$ buys two edges $uv_1,uv_2 \in X_2$, then $d(r,u) \geq 3$. If, in addition, $|T(uv_1)\cup T(uv_2)| > \frac{n}{4}$ then $d(r,u) = 3$. 
\end{observation}

\begin{proof}
By \Cref{7cycle}, $r$ cannot buy edges of $X_2$, therefore $u\not=r$. 
Further, by \Cref{maxn2}, $|T(u)|\leq \frac{n}{2}$. 

We now apply \Cref{mainstrategy1} with the strategy change where $u$ sells $uv_1, uv_2$. Its cost changes by at most
\begin{align*}
&d(r,u)\cdot n -2\alpha -2d(r,u)\cdot |T(u)| + (2d(r,u)+4)\cdot|T(uv_1)\cup T(uv_2)|\\
    &\leq  d(r,u)\cdot n-2\alpha + 4\cdot|T(uv_1)\cup T(uv_2)|\\
    &\leq  d(r,u)\cdot n-2\alpha + 4\frac{n}{2}\\
    &<  d(r,u)\cdot n-2(2n) + 2n \\
    &= (d(r,u)-2)\cdot n
\end{align*}
If $d(r,u) \leq 2$ then this cost change is negative, a contradiction. Thus $d(r,u) \geq 3$.

Now suppose $|T(uv_1)\cup T(uv_2)| > \frac{n}{4}$. Without loss of generality, let $|T(uv_1)| \geq |T(uv_2)|$. Therefore $|T(uv_1)| > \frac{n}{8}$, meaning $T(uv_1) = T(v_1)$. 
As $uv_1 \in X_2$, we know $v_1$ buys an edge in $X_0$ or $X_1\setminus X_0$. Suppose $v_1$ buys $v_1y \in X_0$.  By \Cref{mainstrategy2}, $v$ selling $v_1y$ and buying $v_1r$ changes $v_1$'s cost by
\begin{align*}
    &\leq n - (d(v_1,r)+1)\cdot |T(v_1)| - (d(v_1,r)-1)\cdot|T(u)\setminus T(v_1)|\\
    &\leq n - 2\cdot |T(v_1)| - (d(v,r)-1)\cdot |T(u)| \\
    &< n - 2\frac{n}{8} - (d(v_1,r)-1)\frac{n}{4}\\
    &= n -d(v_1,r)\frac{n}{4}\\
    &= n -(d(u,r)+1)\frac{n}{4}\\
    &\leq n -(3+1)\frac{n}{4}\\
    &=0
\end{align*}
This contradicts the equilibrium condition, therefore
 $v_1$ buys $v_1w \in X_1\setminus X_0$. $w$ buys $wy \in X_0$. By \Cref{mainstrategy2}, $w$ selling $wy$ and buying $wr$ changes $w$'s cost by
\begin{align*}
    &\leq n - (d(w,r)+1)|T(w)|- (d(w,r)-1)|T(v_1)\setminus T(w)| - (d(w,r)-3)|T(u)\setminus T(v_1)|\\
    &\leq n - 2\cdot|T(v_1)| - (d(w,r)-3)\cdot|T(u)| \\
    &< n - 2\frac{n}{8} - (d(w,r)-3)\frac{n}{4}\\
    &= \frac{3}{2}n -d(w,r)\frac{n}{4}\\
    &= \frac{3}{2}n -(d(u,r)+2)\frac{n}{4}\\
    &= n -d(u,r)\frac{n}{4}
\end{align*}
If  $d(r,u) \geq 4$, this leads to a contradiction. Thus  $d(r,u) = 3$, as desired. 

\end{proof}

\section{The Tree Conjecture holds for \texorpdfstring{$\alpha>2n$}{alpha > 2n}}\label{sec:result}

We now have all the resources necessary to prove the tree conjecture holds if $\alpha>2n$.
To do this, we prove two final lemmas which together give the main result.
The first of these lemmas investigates the path to the root $r$ from a vertex $u\in H$
that buys an out-edge.

\begin{lemma}\label{mainlemma1}
 Let $\alpha > 2n$. If $u_0$ buys $u_0v \in X_0$ and $u_0,u_1,...u_{d(u,r)-1},u_{d(u,r)}=r$ is the path from $u_0$ to $r$ in $T$ then at most one of $u_0,u_1,u_2$ has $\deg_H = 2$
\end{lemma}

\begin{proof}
By Observation~\ref{X1}, no vertex buys two edges in $X^+_1$. Thus $u_0$ cannot buy $u_0u_1$ as well as $u_0v$. Therefore $u_1$ buys $u_1u_0$. Similarly, by Observation~\ref{X1}, $u_1$ cannot buy $u_1u_2$ as well as $u_1u_0$.
Therefore $u_2$ buys $u_2u_1$.
    
  Now if $\deg_H(u_0) = 2$ then $S_H(u_0) = T(u_0)$, because $u_1$ is $u_0$'s parent and $v\not\in T(u_0)$, thus $T(u_0)\cap H =\{u_0\}$. Similarly, if $\deg_H(u_1) = 2$ then $S_H(u_1) = T(u_1)\setminus T(u_0)$, and if $\deg_H(u_2) = 2$ then $S_H(u_2) = T(u_2)\setminus T(u_1)$.
    
    Next by \Cref{X2position}, we know that $d(r,u_2) \geq \lfloor\textstyle\frac{|C(u_2)|}{2}\rfloor-2$, 
   that $d(r,u_1) \geq \lfloor\textstyle\frac{|C(u_1)|}{2}\rfloor-1$ and that
 $d(r,u_0) \geq \lfloor\textstyle\frac{|C(u_0)|}{2}\rfloor$.
   
Furthermore, if $\deg_H(u_1) =2$ then, by~\Cref{deg2}, $|S_H(u_1)| \geq |T(u_0)|$.
If $\deg_H(u_2) =2$ then, by~\Cref{directedcycles}, $u_2u_1$ is in a directed cycle, which means $u_3$ buys $u_3u_2$.  Thus, by~\Cref{deg2},  we have $|S_H(u_2)| \geq |T(u_1)|$

    We proceed by case analysis.
    First, assume that $\deg_H(u_2) = \deg_H(u_1) =2$. Then, by~\Cref{deg2}, we have $|S_H(u_1)| \geq \textstyle\frac{\alpha}{2(|C(u_1)|-3)}$ and
    $|S_H(u_2)| \geq |T(u_1)|$. Thus  $|T(u_2)| = |S_H(u_2)| + |T(u_1)| \geq 2|T(u_1)| \geq 2|S_H(u_1)| \geq \textstyle\frac{2\alpha}{2(|C(u_1)|-3)}$. Now, by~\Cref{maxn2}, $|T(u_2)|\leq \textstyle\frac{n}{2}$. Therefore $\textstyle\frac{n}{2} \geq \textstyle\frac{2\alpha}{2(|C(u_1)|-3)}$. Rearranging gives $|C(u_1)| \geq 8$. By Observation~\ref{X2depth}, $d(u_1,r) \geq 3$.
    
   Applying \Cref{mainstrategy2} when $u_1$ sells $u_1u_0$ and buys $u_1r$, the change in cost to $u_1$ is at most 
     \begin{align*}
      &n-|T_{u_\frac{d(u,r)}{2}}| -\textstyle\sum\limits_{l < \textstyle\frac{d(u,r)}{2}}2\cdot|T(u_l)|-(j-1) \cdot\alpha+ \textstyle\textstyle\sum\limits_{k=1}^j (2i_k+d(u,r)+1)\cdot|T(uv_{i_k})|\\
    &=n-|T_{u_\frac{d(u,r)}{2}}| -\textstyle\sum\limits_{l < \textstyle\frac{d(u,r)}{2}}2\cdot|T(u_l)|+ (d(u_1,r)+1+2)\cdot|T(u_1u_0)|\\
    &=n-|T_{u_\frac{d(u,r)}{2}}| -\textstyle\sum\limits_{l < \textstyle\frac{d(u,r)}{2}}2\cdot|T(u_l)|+ (d(u_1,r)+3)\cdot|T(u_0)|\\
    &\leq n-|T_{u_\frac{d(u,r)}{2}}| -\textstyle\sum\limits_{l < \textstyle\frac{d(u,r)}{2}}2\cdot|T(u_l)|+ (2d(u_1,r))\cdot|T(u_0)|\\
    &\leq n-2|T(u_1)|-(d(u_1,r)-1)\cdot|T(u_2)|+ (2d(u_1,r))\cdot|T(u_0)|\\
    &= n-(d(u_1,r)+1)\cdot|T(u_1)|-(d(u_1,r)-1)\cdot(|T(u_2)\setminus|T(u_1)|)+ (2d(u_1,r))\cdot|T(u_0)|\\
    &=n - (d(r,u_1)+1)\cdot |T(u_1)| - (d(u_1,r)-1)\cdot|S_H(u_2)| + (d(u_1,r)+3)\cdot|T(u_0)|\\
    &\leq n - (d(r,u_1)+1)\cdot |T(u_1)| - (d(u_1,r)-1)\cdot|T(u_1)| + (2d(u_1,r))\cdot|T(u_0)|\\
    &= n - 2(d(r,u_1))\cdot |T(u_1)| + 2(d(u_1,r))\cdot|T(u_0)|\\
    &= n - 2(d(r,u_1))\cdot |T(u_1)\setminus T(u_0)|\\
    &= n - 2(d(r,u_1))\cdot |S_H(u_1)|\\
    &\leq n - 2(d(r,u_1))\cdot \textstyle\frac{\alpha}{2(|C(u_1)|-3)}\\
    &\leq n - 2(\lfloor\textstyle\frac{|C(u_1)|}{2}\rfloor-1)\cdot \textstyle\frac{\alpha}{2(|C(u_1)|-3)}\\
    &= n - \textstyle\frac{\alpha}{2}\\
    &< 0
\end{align*}
Therefore it cannot be that $\deg_H(u_2) = \deg_H(u_1) =2$.

Second, suppose $\deg_H(u_0)= \deg_H(u_2)= 2$.
By \Cref{deg2}, we have $|S_H(u_0)| \geq \textstyle\frac{\alpha}{2(|C(u_1)|-3)}$ and $|S_H(u_2)| \geq |T(u_1)|$. 
Because $u_0v \in X_0$, it folows that $T(v) \not\subseteq T(u_0)$. Thus $T(u_0v)=\emptyset$.
Applying \Cref{mainstrategy2} when $u_0$ sells $u_0v$ and buys $u_0r$, the change in cost to $u_0$ is at most
\begin{align*}
      &n-|T_{u_\frac{d(u,r)}{2}}| -\textstyle\sum\limits_{l < \textstyle\frac{d(u,r)}{2}}2\cdot|T(u_l)|-(j-1) \cdot\alpha+ \textstyle\sum\limits_{k=1}^j (2i_k+d(u,r)+1)\cdot|T(uv_{i_k})|\\
    &=n-|T_{u_\frac{d(u,r)}{2}}| -\textstyle\sum\limits_{l < \textstyle\frac{d(u,r)}{2}}2\cdot|T(u_l)|\\
    &\leq n-2|T(u_0)|-2|T(u_1)|-(d(u_0,r)-3)\cdot|T(u_2)|\\
    &= n-2|T(u_0)|-(d(u_0,r)-1)\cdot|T(u_1)|-(d(u_0,r)-3)\cdot|T(u_2)\setminus T(u_1)|\\
    &= n-2|T(u_0)|-(d(u_0,r)-1)\cdot|T(u_1)|-(d(u_0,r)-3)\cdot|S_H(u_2)|\\
    &\leq n-2|T(u_0)|-(d(u_0,r)-1)\cdot|T(u_1)|-(d(u_0,r)-3)\cdot|T(u_1)|\\
    &= n-2|T(u_0)|-(2d(u_0,r)-4)\cdot|T(u_1)|\\
    &\leq n-(2d(u_0,r)-2)|T(u_0)|\\
    &= n-(2d(u_0,r)-2)|S_H(u_0)|\\
    &\leq n - (2d(r,u_0)-2)\cdot\textstyle\frac{\alpha}{2(|C(u_0)|-3)}\\
    &\leq n - (2\lfloor\textstyle\frac{|C(u_0)|}{2}\rfloor-2)\cdot\textstyle\frac{\alpha}{2(|C(u_0)|-3)}\\
    &\leq n - (|C(u_0)|-3)\cdot\textstyle\frac{\alpha}{2(|C(u_0)|-3)}\\
    &= n - \textstyle\frac{\alpha}{2}\\
    &< 0
\end{align*}
Therefore it is not true that $\deg_H(u_2) = \deg_H(u_0) =2$.

Third, suppose $\deg_H(u_0)= \deg_H(u_1)= 2$.
By \Cref{deg2}, we have $|S_H(u_0)| \geq \textstyle\frac{\alpha}{2(|C(u_1)|-3)}$ and
$|S_H(u_1)| \geq |T(u_0)|$.
Because $u_0v \in X_0$, it follows that $T(v) \not\subseteq T(u_0)$. Thus $T(u_0v)=\emptyset$.
Applying \Cref{mainstrategy2} when $u_0$ sells $u_0v$ and buys $u_0r$, the change in cost to $u_0$ is
\begin{align*}
      &\leq n-|T_{u_\frac{d(u,r)}{2}}| -\textstyle\sum\limits_{l < \textstyle\frac{d(u,r)}{2}}2\cdot|T(u_l)|-(j-1) \cdot\alpha+ \textstyle\sum\limits_{k=1}^j (2i_k+d(u,r)+1)\cdot|T(uv_{i_k})|\\
    &=n-|T_{u_\frac{d(u,r)}{2}}| -\textstyle\sum\limits_{l < \textstyle\frac{d(u,r)}{2}}2\cdot|T(u_l)|\\
    &\leq n-2|T(u_0)|-(d(u_0,r)-1)\cdot|T(u_1)|\\
    &= n-(d(u_0,r)+1)\cdot|T(u_0)|-(d(u_0,r)-1)\cdot|T(u_1)\setminus T(u_0)|\\
    &= n-(d(u_0,r)+1)\cdot|T(u_0)|-(d(u_0,r)-1)\cdot|S_H(u_1)|\\
    &\leq n-(d(u_0,r)+1)\cdot|T(u_0)|-(d(u_0,r)-1)\cdot|T(u_0)|\\
    &\leq n-(2d(u_0,r))\cdot|T(u_0)|\\
    &= n-(2d(u_0,r))\cdot|S_H(u_0)|\\
    &\leq n - (2d(r,u_0))\cdot\textstyle\frac{\alpha}{2(|C(u_0)|-3)}\\
    &\leq n - (2\lfloor\textstyle\frac{|C(u_0)|}{2}\rfloor)\cdot\textstyle\frac{\alpha}{2(|C(u_0)|-3)}\\
    &\leq n - (|C(u_0)|-3)\cdot\textstyle\frac{\alpha}{2(|C(u_0)|-3)}\\
    &= n - \textstyle\frac{\alpha}{2}\\
    &< 0
\end{align*}
  
    Therefore it is not true that $\deg_H(u_1) = \deg_H(u_0) =2$.
    Thus, at most one of $u_0,u_1,u_2$ has $\deg_H = 2$, as desired.
\end{proof}

We now apply a counting argument to upper bound the number of vertices in $H$ that buy two edges in $X_2$.
\begin{lemma}\label{mainlemma2}
    The number of vertices $v \in H$ that buy two edges in $X_2$ is $< \deg_H(r)$.
\end{lemma}

\begin{proof}
The only vertices that can buy two edges $e_1,e_2$ in $X_2$ are those with $|T(e_1)\cup T(e_2)|>\textstyle\frac{n}{4}$ by Observation~\ref{X2}. 
All of these vertices are distance 3 from $r$, by Observation~\ref{X2depth}. The common ancestor in $T$ of any pair of these vertices must be $r$, otherwise there is a vertex $u\not= r$  with $|T(u)| >\textstyle\frac{n}{2}$, contradicting \Cref{maxn2}. Furthermore, $r$ cannot buy any edge $e$ such that $T(e)$ contains such a vertex $v$, otherwise $r$ could sell $e$ and buy $rv$, reducing its distance to $T(v)$ by 2 and increasing its distance to $T(e)\setminus T(v)$ by $\leq 2$. This gives a net increase in cost of $\leq 2|T(v)|-2|T(e)\setminus T(v)| \leq 4|T(v)| -  2|T(e)| < 4\textstyle\frac{n}{4}- 2\textstyle\frac{n}{2} = 0$, a contradiction. Therefore the number of vertices which buy two edges in $X_2$ is $\leq deg_H^-(r)$. By \Cref{directedcycles}, we must have $deg_H^+(r)\geq 1$, implying $deg_H^-(r) < deg_H(r)$ and completing the proof.

\end{proof}

Putting everything together we can now prove the main result.
\begin{proof}[of~\Cref{main}]
    Here we use combine all to prove that if $G$ is a Nash equilibrium graph for the network creation game $(n,\alpha)$ and $\alpha > 2n$, $G$ is a tree.

    This is a proof by contradiction based on the assumption that there exists a biconnected component $H$ in $G$ containing $n_H$ vertices. 
    Within $H$, the sum of degrees of all vertices in $H$ equals $2(n_H-1)+ 2|X_0|$. That is twice the number of edges in the spanning tree on $H$ induced by $T$ plus twice the number of out-edges in $H$. 

    Let's instead count the degrees in the following way: 
    \begin{align*}
    &\textstyle\sum\limits_{i=2}^{n_H}\, \textstyle\sum\limits_{v\in H:\deg_H(v)=i}i\\
    &= 2n_H + \textstyle\sum\limits_{i=3}^{n_H}\, \textstyle\sum\limits_{v\in H:\deg_H(v)=i}i-2\\
    &= 2n_H + \deg_H(r)-2 + \textstyle\sum\limits_{i=3}^{n_H}\,\textstyle\sum\limits_{v\in H-r:\deg_H(v)=i}i-2\\
    &\geq 2n_H + \deg_H(r)-2 + 2|X_0|-\textstyle\sum\limits_{i=2}^{|X_2|}\,\textstyle\sum\limits_{v \in H: \text{$v$ buys $i$ edges in $X_2$}}i-1\quad{\text{[by \Cref{mainlemma1}]}}\\
    &\geq 2n_H + \deg_H(r)-2 + 2|X_0|-\textstyle\sum\limits_{v \in H: \text{$v$ buys $2$ edges in $X_2$}}1\qquad{\text{[by Observation~\ref{X2}]}}\\
    &> 2n_H + \deg_H(r)-2 + 2|X_0| - \deg_H(r)\qquad{\text{[by \Cref{mainlemma2}]}}\\
    &= 2(n_H-1)+2|X_0|
    \end{align*}  
This is a contradiction, implying that $H$ does not exist for $\alpha > 2n$.
    
\end{proof}

\section{Conclusion}\label{sec:conclusion}

In this paper we proved the revised tree conjecture holds for $\alpha >2n$.
Moreover, we have reached a natural limit in the quest to settle this conjecture.
Specifically, for $\alpha\in [n,2n]$, the range in which the conjecture is unsettled, 
we are no longer certain that directed cycles must be present in non-tree equilibrium graphs. 
To close his remaining gap, we believe min-cycles still have an important role to play but, to allow 
their usage, more precise analyses will be necessary. A potentially useful intermediate step would 
be to determine conditions that allow for an undirected min-cycle.

\newpage
\bibliography{main}

\end{document}